\documentclass[11pt]{article}

\usepackage{latexsym}
\usepackage{amssymb}

\newtheorem{thm}{Theorem}

\newtheorem{prop}{Proposition}

\newtheorem{dfn}{Definition}

\begin{document}

\begin{center}
{\Large\bf
The matrix Stieltjes moment problem: a description of all solutions.}
\end{center}
\begin{center}
{\bf S.M. Zagorodnyuk}
\end{center}

\section{Introduction.}
The matrix Stieltjes moment problem consists of
finding a left-continuous non-decreasing matrix function $M(x) = ( m_{k,l}(x) )_{k,l=0}^{N-1}$
on $\mathbb{R}_+ = [0,+\infty)$, $M(0)=0$, such that
\begin{equation}
\label{f1_1}
\int_{\mathbb{R}_+} x^n dM(x) = S_n,\qquad n\in\mathbb{Z}_+,
\end{equation}
where $\{ S_n \}_{n=0}^\infty$ is a given sequence of Hermitian $(N\times N)$ complex matrices, $N\in\mathbb{N}$.
This problem is said to be determinate, if there exists a unique solution and indeterminate in the opposite case.

In the scalar ($N=1$) indeterminate case the Stieltjes moment problem was solved by M.G.~Krein
(see~\cite{c_1000_Kr_st},\cite{c_2000_KrN}), while in the scalar degenerate case the problem was solved
by F.R.~Gantmacher in~\cite[Chapter XVI]{Cit_3000_Gantmacher}.

The operator (and, in particular, the matrix) Stieltjes moment problem was introduced by M.G.~Krein and
M.A.~Krasnoselskiy in~\cite{c_4000_Kr_Kras}. They obtained the necessary and sufficient  conditions of
solvability for this problem.

\noindent
Let us introduce the following matrices
\begin{equation}
\label{f1_2}
\Gamma_n = (S_{i+j})_{i,j=0}^n = \left(
\begin{array}{cccc} S_0 & S_1 & \ldots & S_n\\
S_1 & S_2 & \ldots & S_{n+1}\\
\vdots & \vdots & \ddots & \vdots\\
S_n & S_{n+1} & \ldots & S_{2n}\end{array}
\right),
\end{equation}
\begin{equation}
\label{f1_3}
\widetilde\Gamma_n = (S_{i+j+1})_{i,j=0}^n = \left(
\begin{array}{cccc} S_1 & S_2 & \ldots & S_{n+1}\\
S_2 & S_3 & \ldots & S_{n+2}\\
\vdots & \vdots & \ddots & \vdots\\
S_{n+1} & S_{n+2} & \ldots & S_{2n+1}\end{array}
\right),\qquad n\in\mathbb{Z}_+.
\end{equation}
The moment problem~(\ref{f1_1}) has a solution if and only if
\begin{equation}
\label{f1_4}
\Gamma_n \geq 0,\quad \widetilde\Gamma_n \geq 0,\qquad n\in\mathbb{Z}_+.
\end{equation}
In~2004, Yu.M.~Dyukarev performed a deep investigation of the moment problem~(\ref{f1_1}) in
the case when
\begin{equation}
\label{f1_5}
\Gamma_n > 0,\quad \widetilde\Gamma_n > 0,\qquad n\in\mathbb{Z}_+,
\end{equation}
and some limit matrix intervals (which he called the limit Weyl intervals) are non-degenerate, see~\cite{c_5000_D}.
He obtained a parameterization of all solutions of the moment problem in this case.

\noindent
Our aim here is to obtain a description of all solutions of the moment problem~(\ref{f1_1}) in the
general case. No conditions besides the solvability (i.e. conditions~(\ref{f1_4})) will be assumed.
We shall apply an operator approach which was used in~\cite{c_6000_Z} and Krein's formula for
the generalized $\Pi$-resolvents of non-negative Hermitian operators~\cite{c_7000_Kr},\cite{c_8000_Kr_Ovch}. We shall
use Krein's formula in the form which was proposed by V.A.~Derkach and
M.M.~Malamud in~\cite{c_9000_D_M}. We should also notice that these authors
presented a detailed proof of Krein's formula.

\noindent
{\bf Notations.}  As usual, we denote by $\mathbb{R}, \mathbb{C}, \mathbb{N}, \mathbb{Z}, \mathbb{Z}_+$
the sets of real numbers, complex numbers, positive integers, integers and non-negative integers,
respectively; $\mathbb{R}_+ = [0,+\infty)$, $\mathbb{C}_+ = \{ z\in \mathbb{C}:\ \mathop{\rm Im}\nolimits z>0 \}$.
The space of $n$-dimensional complex vectors
$a = (a_0,a_1,\ldots,a_{n-1})$, will be
denoted by $\mathbb{C}^n$, $n\in \mathbb{N}$.
If $a\in \mathbb{C}^n$ then $a^*$ means the complex conjugate vector.
By $\mathbb{P}$ we denote the set of all complex polynomials.

\noindent
Let $M(x)$ be a left-continuous non-decreasing matrix function $M(x) = ( m_{k,l}(x) )_{k,l=0}^{N-1}$
on $\mathbb{R_+}$, $M(0)=0$, and $\tau_M (x) := \sum_{k=0}^{N-1} m_{k,k} (x)$;
$\Psi(x) = ( dm_{k,l}/ d\tau_M )_{k,l=0}^{N-1}$ (the Radon-Nikodym derivative).
We denote by $L^2(M)$ a set (of classes of equivalence)
of vector functions $f: \mathbb{R}\rightarrow \mathbb{C}^N$, $f = (f_0,f_1,\ldots,f_{N-1})$,
such that (see, e.g.,~\cite{c_10000_M_M})
$$ \| f \|^2_{L^2(M)} := \int_\mathbb{R}  f(x) \Psi(x) f^*(x) d\tau_M (x) < \infty. $$
The space $L^2(M)$ is a Hilbert space with the scalar product
$$ ( f,g )_{L^2(M)} := \int_\mathbb{R}  f(x) \Psi(x) g^*(x) d\tau_M (x),\qquad f,g\in L^2(M). $$

For a separable Hilbert space $H$ we denote by $(\cdot,\cdot)_H$ and $\| \cdot \|_H$ the scalar
product and the norm in $H$, respectively. The indices may be omitted in obvious cases.
By $E_H$ we denote the identity operator in $H$, i.e. $E_H x = x$, $x\in H$.

\noindent
For a linear operator $A$ in $H$ we denote by $D(A)$ its domain, by $R(A)$ its range, and by
$\mathop{\rm ker}\nolimits A$ its kernel. By $A^*$ we denote its adjoint if it exists.
By $\rho(A)$ we denote the resolvent set of $A$; $N_z = \mathop{\rm ker}\nolimits(A^* - zE_H)$.
If $A$ is bounded, then $\| A \|$ stands for its operator norm.
For a set of elements $\{ x_n \}_{n\in T}$ in $H$, we
denote by $\mathop{\rm Lin}\nolimits\{ x_n \}_{n\in T}$ and $\mathop{\rm span}\nolimits\{ x_n \}_{n\in T}$ the linear
span and the closed
linear span (in the norm of $H$), respectively. Here $T$ is an arbitrary set of indices.
For a set $M\subseteq H$ we denote by $\overline{M}$ the closure of $M$ with respect to the norm of $H$.

\noindent
If $H_1$ is a subspace of $H$, by $P_{H_1} = P_{H_1}^{H}$ we denote the operator of the orthogonal projection on $H_1$
in $H$. If $\mathcal{H}$ is another Hilbert space, by $[H,\mathcal{H}]$ we denote the space of all
bounded operators from $H$ into $\mathcal{H}$; $[H]:= [H,H]$.
$\mathfrak{C}(H)$ is the set of closed linear operators $A$ such that $\overline{D(A)}=H$.

\section{The matrix Stieltjes moment problem: the solvability.}
Consider the matrix Stieltjes moment problem~(\ref{f1_1}).
Let us check that conditions~(\ref{f1_4}) are necessary for the solvability
of the problem~(\ref{f1_1}). In  fact, suppose that the moment problem has
a solution $M(x)$. Choose an arbitrary function
$a(x) = (a_0(x),a_1(x),...,a_{N-1}(x))$, where
$$ a_j(x) = \sum_{k=0}^n \alpha_{j,k} x^k,\quad \alpha_{j,k}\in \mathbb{C},\ n\in \mathbb{Z}_+. $$
This function belongs to $L^2(M)$ and
$$ 0 \leq \int_{\mathbb{R}_+} a(x) dM(x) a^*(x) =
\sum_{k,l=0}^n \int_{\mathbb{R}_+}(\alpha_{0,k},\alpha_{1,k},...,\alpha_{N-1,k}) x^{k+l} dM(x) $$
$$* (\alpha_{0,l},\alpha_{1,l},...,\alpha_{N-1,l})^* =
\sum_{k,l=0}^n (\alpha_{0,k},\alpha_{1,k},...,\alpha_{N-1,k}) S_{k+l} $$
$$* (\alpha_{0,l},\alpha_{1,l},...,\alpha_{N-1,l})^* =
A\Gamma_n A^*, $$
where $A = (\alpha_{0,0},\alpha_{1,0},...,\alpha_{N-1,0},\alpha_{0,1},\alpha_{1,1},...,
\alpha_{N-1,1},...,\alpha_{0,n},\alpha_{1,n},...,\alpha_{N-1,n})$,
and we have used the rules for the multiplication of block matrices.
In a similar manner we get
$$ 0 \leq \int_{\mathbb{R}_+} a(x) x dM(x) a^*(x) = A\widetilde\Gamma_n A^*, $$
and therefore conditions~(\ref{f1_4}) hold.

\noindent
On the other hand, let the moment problem~(\ref{f1_1}) be given and suppose  that
conditions~(\ref{f1_4}) are true.
For the prescribed moments
$$ S_j = (s_{j;k,l})_{k,l=0}^{N-1},\quad s_{j;k,l}\in \mathbb{C},\qquad j\in \mathbb{Z}_+, $$
we consider the following block matrices
\begin{equation}
\label{f2_1}
\Gamma = (S_{i+j})_{i,j=0}^\infty =
\left(
\begin{array}{cccc}
S_{0} & S_{1} & S_2 & \ldots \\
S_{1} & S_{2} & S_3 & \ldots \\
S_{2} & S_{3} & S_4 & \ldots \\
\vdots & \vdots & \vdots & \ddots \end{array}\right),
\end{equation}
\begin{equation}
\label{f2_1_1}
\widetilde\Gamma = (S_{i+j+1})_{i,j=0}^\infty =
\left(
\begin{array}{cccc}
S_{1} & S_{2} & S_3 & \ldots \\
S_{2} & S_{3} & S_4 & \ldots \\
S_{3} & S_{4} & S_5 & \ldots \\
\vdots & \vdots & \vdots & \ddots \end{array}\right).
\end{equation}
The matrix $\Gamma$ can be viewed as a scalar semi-infinite matrix
\begin{equation}
\label{f2_2}
\Gamma = (\gamma_{n,m})_{n,m=0}^\infty,\qquad \gamma_{n,m}\in \mathbb{C}.
\end{equation}
Notice that
\begin{equation}
\label{f2_3}
\gamma_{rN+j,tN+n} = s_{r+t;j,n},\qquad r,t\in \mathbb{Z_+},\ 0\leq j,n\leq N-1.
\end{equation}
The matrix $\widetilde\Gamma$ can be also viewed as a scalar semi-infinite matrix
\begin{equation}
\label{f2_3_1}
\widetilde\Gamma = (\widetilde\gamma_{n,m})_{n,m=0}^\infty =
(\gamma_{n+N,m})_{n,m=0}^\infty.
\end{equation}
The conditions in~(\ref{f1_4}) imply that
\begin{equation}
\label{f2_4}
(\gamma_{k,l})_{k,l=0}^r \geq 0,\qquad r\in \mathbb{Z}_+;
\end{equation}
\begin{equation}
\label{f2_4_1}
(\gamma_{k+N,l})_{k,l=0}^r \geq 0,\qquad r\in \mathbb{Z}_+.
\end{equation}
We shall use the following important fact (e.g.,~\cite[Supplement 1]{c_11000_AG}):
\begin{thm}
\label{t2_1}
Let $\Gamma = (\gamma_{n,m})_{n,m=0}^\infty$, $\gamma_{n,m}\in \mathbb{C}$, be a semi-infinite
complex matrix such that condition~(\ref{f2_4}) holds.
Then there exist a separable Hilbert space $H$ with a scalar product $(\cdot,\cdot)_H$ and
a sequence $\{ x_n \}_{n=0}^\infty$ in $H$, such that
\begin{equation}
\label{f2_5}
\gamma_{n,m} = (x_n,x_m)_H,\qquad n,m\in \mathbb{Z}_+,
\end{equation}
and $\mathop{\rm span}\nolimits\{ x_n \}_{n=0}^\infty = H$.
\end{thm}
{\bf Proof. }
Consider an arbitrary infinite-dimensional linear vector space $V$. For example, we can choose the linear space
of all complex sequences $(u_n)_{n\in \mathbb{Z}_+}$, $u_n\in \mathbb{C}$.
Let $X = \{ x_n \}_{n=0}^\infty$ be an arbitrary infinite sequence of linear independent elements
in $V$. Let $L = \mathop{\rm Lin}\nolimits\{ x_n \}_{n\in\mathbb{Z}_+}$ be the linear span of elements of $X$. Introduce the following functional:
\begin{equation}
\label{f2_6}
[x,y] = \sum_{n,m=0}^\infty \gamma_{n,m} a_n\overline{b_m},
\end{equation}
for $x,y\in L$,
$$ x=\sum_{n=0}^\infty a_n x_n,\quad y=\sum_{m=0}^\infty b_m x_m,\quad a_n,b_m\in\mathbb{C}. $$
Here and in what follows we assume that for elements of linear spans all but a finite number of coefficients
are zero.
The space $V$ with $[\cdot,\cdot]$ will be a quasi-Hilbert space. Factorizing and making the completion
we obtain the  required space $H$ (see~\cite{c_12000_Ber}).
$\Box$

From~(\ref{f2_3}) it follows that
\begin{equation}
\label{f2_7}
\gamma_{a+N,b} = \gamma_{a,b+N},\qquad a,b\in\mathbb{Z}_+.
\end{equation}
In fact, if $a=rN+j$, $b=tN+n$, $0\leq j,n \leq N-1$,  $r,t\in\mathbb{Z}_+$, we can write
$$ \gamma_{a+N,b} = \gamma_{(r+1)N+j,tN+n} = s_{r+t+1;j,n} = \gamma_{rN+j,(t+1)N+n} = \gamma_{a,b+N}. $$
By Theorem~\ref{t2_1} there exist a Hilbert space $H$ and a sequence $\{ x_n \}_{n=0}^\infty$ in $H$,
such that $\mathop{\rm span}\nolimits\{ x_n \}_{n=0}^\infty = H$, and
\begin{equation}
\label{f2_8}
(x_n,x_m)_H = \gamma_{n,m},\qquad n,m\in\mathbb{Z}_+.
\end{equation}
Set $L := \mathop{\rm Lin}\nolimits\{ x_n \}_{n=0}^\infty$.
Notice that elements $\{ x_n \}$ are {\it not necessarily linearly independent}. Thus, for an
arbitrary $x\in L$ there can exist different representations:
\begin{equation}
\label{f2_9}
x = \sum_{k=0}^\infty \alpha_k x_k,\quad \alpha_k\in \mathbb{C},
\end{equation}
\begin{equation}
\label{f2_10}
x = \sum_{k=0}^\infty \beta_k x_k,\quad \beta_k\in \mathbb{C}.
\end{equation}
(Here all but a finite number of coefficients $\alpha_k$, $\beta_k$ are zero).
Using~(\ref{f2_7}),(\ref{f2_8}) we can write
$$ \left( \sum_{k=0}^\infty \alpha_k x_{k+N}, x_l \right) = \sum_{k=0}^\infty \alpha_k ( x_{k+N}, x_l ) =
\sum_{k=0}^\infty \alpha_k \gamma_{k+N,l} = \sum_{k=0}^\infty \alpha_k \gamma_{k,l+N} $$
$$ = \sum_{k=0}^\infty \alpha_k ( x_{k}, x_{l+N} ) =
\left( \sum_{k=0}^\infty \alpha_k x_{k}, x_{l+N} \right) = (x,x_{l+N}),\qquad l\in\mathbb{Z}_+. $$
In a similar manner we obtain that
$$ \left( \sum_{k=0}^\infty \beta_k x_{k+N}, x_l \right) = (x,x_{l+N}),\qquad l\in\mathbb{Z}_+, $$
and therefore
$$ \left( \sum_{k=0}^\infty \alpha_k x_{k+N}, x_l \right) =
\left( \sum_{k=0}^\infty \beta_k x_{k+N}, x_l \right),\qquad l\in\mathbb{Z}_+. $$
Since $\overline{L} = H$, we obtain that
\begin{equation}
\label{f2_11}
\sum_{k=0}^\infty \alpha_k x_{k+N} = \sum_{k=0}^\infty \beta_k x_{k+N}.
\end{equation}
Let us introduce the following operator:
\begin{equation}
\label{f2_12}
A x = \sum_{k=0}^\infty \alpha_k x_{k+N},\qquad x\in L,\ x = \sum_{k=0}^\infty \alpha_k x_{k}.
\end{equation}
Relations~(\ref{f2_9}),(\ref{f2_10}) and~(\ref{f2_11}) show that this definition does not
depend on the choice of a representation for $x\in L$. Thus, this definition is correct.
In particular, we have
\begin{equation}
\label{f2_13}
A x_k =  x_{k+N},\qquad k\in\mathbb{Z}_+.
\end{equation}
Choose arbitrary $x,y\in L$, $x = \sum_{k=0}^\infty \alpha_k x_{k}$, $y = \sum_{n=0}^\infty \gamma_n x_{n}$,
and write
$$ (Ax,y) = \left( \sum_{k=0}^\infty \alpha_k x_{k+N},\sum_{n=0}^\infty \gamma_n x_{n} \right) =
\sum_{k,n=0}^\infty \alpha_k \overline{\gamma_n} (x_{k+N},x_n) $$
$$ = \sum_{k,n=0}^\infty \alpha_k \overline{\gamma_n} (x_{k},x_{n+N}) =
\left( \sum_{k=0}^\infty \alpha_k x_{k},\sum_{n=0}^\infty \gamma_n x_{n+N} \right) =
(x,Ay). $$
By relation~(\ref{f2_4_1}) we get
$$ (Ax,x) = \left( \sum_{k=0}^\infty \alpha_k x_{k+N},\sum_{n=0}^\infty \alpha_n x_{n} \right) =
\sum_{k,n=0}^\infty \alpha_k \overline{\alpha_n} (x_{k+N},x_n) $$
$$ = \sum_{k,n=0}^\infty \alpha_k \overline{\alpha_n} \gamma_{k+N,n} \geq 0, $$
Thus, the operator $A$ is a linear non-negative Hermitian operator in $H$ with the domain $D(A)=L$.
Such an operator has a non-negative self-adjoint extension~\cite[Theorem 7, p.450]{c_13000_Kr}.
Let $\widetilde A\supseteq A$ be an arbitrary non-negative self-adjoint extension of $A$ in a Hilbert space
$\widetilde H\supseteq H$, and $\{ \widetilde E_\lambda \}_{\lambda\in \mathbb{R}_+}$ be its left-continuous orthogonal
resolution of unity.
Choose an arbitrary $a\in \mathbb{Z}_+$, $a=rN + j$, $r\in \mathbb{Z}_+$, $0\leq j\leq N-1$. Notice that
$$ x_a = x_{rN+j} = A x_{(r-1)N+j} = ... = A^r x_j. $$
Using~(\ref{f2_3}),(\ref{f2_8}) we can write
$$ s_{r+t;j,n} = \gamma_{rN+j,tN+n} = ( x_{rN+j},x_{tN+n} )_H = (A^r x_j, A^t x_n)_H $$
$$ = ( \widetilde A^r x_j, \widetilde A^t x_n)_{\widetilde H} =
\left( \int_{\mathbb{R}_+} \lambda^r d\widetilde E_\lambda x_j, \int_{\mathbb{R}_+} \lambda^t d\widetilde E_\lambda x_n
\right)_{\widetilde H} $$
$$ = \int_{\mathbb{R}_+} \lambda^{r+t} d (\widetilde E_\lambda x_j, x_n)_{\widetilde H} =
\int_{\mathbb{R}_+} \lambda^{r+t} d \left( P^{\widetilde H}_H \widetilde E_\lambda x_j, x_n \right)_{H}. $$
Let us write the last relation in a matrix form:
\begin{equation}
\label{f2_14}
S_{r+t} = \int_{\mathbb{R}_+} \lambda^{r+t} d \widetilde M(\lambda),\qquad r,t\in\mathbb{Z}_+,
\end{equation}
where
\begin{equation}
\label{f2_15}
\widetilde M(\lambda) := \left( \left( P^{\widetilde H}_H \widetilde E_\lambda x_j,
x_n \right)_{H} \right)_{j,n=0}^{N-1}.
\end{equation}
If we set $t=0$ in relation~(\ref{f2_14}), we obtain that the matrix function $\widetilde M(\lambda)$ is
a solution of the matrix Stieltjes moment  problem~(\ref{f1_1}). In fact, from the properties of the
orthogonal resolution of unity it easily follows that $\widetilde M (\lambda)$ is left-continuous non-decreasing and
$\widetilde M(0) = 0$.

Thus, we obtained another proof of the solvability criterion for the matrix Stieltjes moment  problem~(\ref{f1_1}):
\begin{thm}
\label{t2_2}
Let a matrix Stieltjes moment problem~(\ref{f1_1}) be given. This problem has a solution if and only if
conditions~(\ref{f1_4}) hold true.
\end{thm}

\section{A description of solutions.}
Let $B$ be an arbitrary non-negative Hermitian operator in a Hilbert space $\mathcal{H}$.
Choose an arbitrary non-negative self-adjoint extension $\widehat B$ of $B$ in a Hilbert space
$\widehat{\mathcal{H}} \supseteq \mathcal{H}$.
Let $R_z(\widehat B)$ be the resolvent of $\widehat B$ and $\{ \widehat E_\lambda\}_{\lambda\in \mathbb{R}_+}$
be the orthogonal left-continuous resolution of unity of $\widehat B$. Recall that the operator-valued function
$\mathbf R_z = P_{ \mathcal{H} }^{ \widehat{\mathcal{H}} } R_z(\widehat B)$ is called {\bf a generalized
$\Pi$-resolvent of $B$}, $z\in\mathbb{C}\backslash\mathbb{R}$~\cite{c_8000_Kr_Ovch}.
If $\widehat{\mathcal{H}} = \mathcal{H}$ then $R_z(\widehat B)$ is called {\bf a canonical $\Pi$-resolvent}.
The function
$\mathbf E_\lambda = P_{\mathcal{H}}^{\widehat{\mathcal{H}}} \widehat E_\lambda$, $\lambda\in\mathbb{R}$, we call
a {\bf $\Pi$-spectral
function} of a non-negative Hermitian operator $B$.
There exists a one-to-one correspondence between generalized $\Pi$-resolvents and $\Pi$-spectral functions
established by the following relation (\cite{c_11000_AG}):
\begin{equation}
\label{f3_1}
(\mathbf R_z f,g)_{\mathcal{H}} = \int_{\mathbb{R}_+} \frac{1}{\lambda - z}
d( \mathbf E_\lambda f,g)_{\mathcal{H}},\qquad f,g\in \mathcal{H},\
z\in \mathbb{C}\backslash \mathbb{R}.
\end{equation}
Denote the set of all generalized $\Pi$-resolvents of $B$ by $\Omega^0(-\infty,0)=\Omega^0(-\infty,0)(B)$.

Let a moment problem~(\ref{f1_1}) be given and conditions~(\ref{f1_4}) hold.
Consider the operator $A$ defined as in~(\ref{f2_12}).
Formula~(\ref{f2_15}) shows that $\Pi$-spectral functions of the operator $A$ produce
solutions of the matrix Stieltjes moment problem~(\ref{f1_1}).
Let us show that an arbitrary solution of~(\ref{f1_1}) can be produced in this way.

\noindent
Choose an arbitrary solution $\widehat M(x) = ( \widehat m_{k,l}(x) )_{k,l=0}^{N-1}$ of
the matrix Stieltjes moment problem~(\ref{f1_1}). Consider the space $L^2(\widehat M)$ and
let $Q$ be the operator of multiplication by an independent variable in $L^2(\widehat M)$.
The operator $Q$ is self-adjoint and its resolution of unity is given by (see~\cite{c_10000_M_M})
\begin{equation}
\label{f3_2}
E_b - E_a = E([a,b)): h(x) \rightarrow \chi_{[a,b)}(x) h(x),
\end{equation}
where $\chi_{[a,b)}(x)$ is the characteristic function of an interval $[a,b)$, $0 \leq a<b\leq +\infty$.
Set
$$ \vec e_k = (e_{k,0},e_{k,1},\ldots,e_{k,N-1}),\quad e_{k,j}=\delta_{k,j},\qquad 0\leq j\leq N-1, $$
where $k=0,1,\ldots N-1$.
A set of (classes of equivalence of) functions $f\in L^2(\widehat M)$ such that
(the corresponding class includes) $f=(f_0,f_1,\ldots, f_{N-1})$, $f\in \mathbb{P}$, we denote
by $\mathbb{P}^2(\widehat M)$. It is said to be a set of vector polynomials in $L^2(\widehat M)$.
Set $L^2_0(\widehat M) := \overline{ \mathbb{P}^2(\widehat M) }$.

For an arbitrary (representative) $f\in \mathbb{P}^2(\widehat M)$ there
exists a unique representation of the following form:
\begin{equation}
\label{f3_3}
f(x) = \sum_{k=0}^{N-1}
\sum_{j=0}^\infty \alpha_{k,j} x^j \vec e_k,\quad \alpha_{k,j}\in \mathbb{C}.
\end{equation}
Here the sum is assumed to be finite.
Let $g\in \mathbb{P}^2(\widehat M)$ have a representation
\begin{equation}
\label{f3_4}
g(x) = \sum_{l=0}^{N-1} \sum_{r=0}^\infty \beta_{l,r} x^r \vec e_l,\quad \beta_{l,r}\in \mathbb{C}.
\end{equation}
Then we can write
$$ (f,g)_{L^2(\widehat M)} = \sum_{k,l=0}^{N-1} \sum_{j,r=0}^\infty \alpha_{k,j}\overline{\beta_{l,r}}
\int_\mathbb{R} x^{j+r} \vec e_k d\widehat M(x) \vec e_l^*  $$
\begin{equation}
\label{f3_5}
= \sum_{k,l=0}^{N-1}
\sum_{j,r=0}^\infty \alpha_{k,j}\overline{\beta_{l,r}}
\int_\mathbb{R} x^{j+r} d\widehat m_{k,l}(x)
= \sum_{k,l=0}^{N-1} \sum_{j,r=0}^\infty \alpha_{k,j}\overline{\beta_{l,r}}
s_{j+r;k,l}.
\end{equation}
On the other hand, we can write
$$ \left( \sum_{j=0}^\infty \sum_{k=0}^{N-1} \alpha_{k,j} x_{jN+k},
\sum_{r=0}^\infty \sum_{l=0}^{N-1} \beta_{l,r} x_{rN+l} \right)_H =
\sum_{k,l=0}^{N-1} \sum_{j,r=0}^\infty \alpha_{k,j}\overline{\beta_{l,r}}
(x_{jN+k}, x_{rN+l})_H  $$
\begin{equation}
\label{f3_6}
= \sum_{k,l=0}^{N-1} \sum_{j,r=0}^\infty \alpha_{k,j}\overline{\beta_{l,r}}
\gamma_{jN+k,rN+l}
= \sum_{k,l=0}^{N-1} \sum_{j,r=0}^\infty \alpha_{k,j}\overline{\beta_{l,r}}
s_{j+r;k,l}.
\end{equation}
From relations~(\ref{f3_5}),(\ref{f3_6}) it follows that
\begin{equation}
\label{f3_7}
(f,g)_{L^2(\widehat M)} = \left( \sum_{j=0}^\infty \sum_{k=0}^{N-1} \alpha_{k,j} x_{jN+k},
\sum_{r=0}^\infty \sum_{l=0}^{N-1} \beta_{l,r} x_{rN+l} \right)_H.
\end{equation}
Let us introduce the following operator:
\begin{equation}
\label{f3_8}
Vf = \sum_{j=0}^\infty \sum_{k=0}^{N-1} \alpha_{k,j} x_{jN+k},
\end{equation}
for $f(x)\in \mathbb{P}^2(\widehat M)$, $f(x) = \sum_{k=0}^{N-1} \sum_{j=0}^\infty \alpha_{k,j} x^j \vec e_k$,
$\alpha_{k,j}\in \mathbb{C}$.
Let us show that this definition is correct. In fact,
if vector polynomials $f$, $g$ have representations~(\ref{f3_3}),(\ref{f3_4}), and
$\| f-g \|_{L^2(\widehat M)} = 0$, then
from~(\ref{f3_7}) it follows that $V(f-g)=0$.
Thus, $V$ is a correctly defined operator from $\mathbb{P}^2(\widehat M)$ into $H$.

Relation~(\ref{f3_7}) shows that $V$ is an isometric transformation from $\mathbb{P}^2(\widehat M)$ onto $L$.
By continuity we extend it to an isometric transformation from $L^2_0(\widehat M)$ onto $H$.
In particular, we note that
\begin{equation}
\label{f3_9}
V x^j \vec e_k = x_{jN+k},\qquad j\in \mathbb{Z}_+;\quad 0\leq k\leq N-1.
\end{equation}
Set $L^2_1 (\widehat M) := L^2(\widehat M)\ominus L^2_0 (\widehat M)$, and
$U := V\oplus E_{L^2_1 (\widehat M)}$. The operator $U$ is
an isometric transformation from $L^2(\widehat M)$ onto $H\oplus L^2_1 (\widehat M)=:\widehat H$.
Set
$$ \widehat A := UQU^{-1}. $$
The operator $\widehat A$ is a non-negative self-adjoint operator in $\widehat H$.
Let $\{ \widehat E_\lambda \}_{\lambda\in\mathbb{R}_+}$
be its left-continuous orthogonal resolution of unity.
Notice that
$$ UQU^{-1} x_{jN+k} = VQV^{-1} x_{jN+k} = VQ x^j \vec e_k = V x^{j+1} \vec e_k =
x_{(j+1)N+k} $$
$$ = x_{jN+k+N} = Ax_{jN+k},\qquad j\in\mathbb{Z}_+;\quad 0\leq k\leq N-1. $$
By linearity we get
$$ UQU^{-1} x = Ax,\qquad x\in L = D(A), $$
and therefore $\widehat A\supseteq A$.
Choose an arbitrary $z\in\mathbb{C}\backslash\mathbb{R}$ and write
$$ \int_{\mathbb{R}_+} \frac{1}{\lambda - z} d( \widehat E_\lambda x_k, x_j)_{\widehat H} =
\left( \int_{\mathbb{R}_+} \frac{1}{\lambda - z} d\widehat E_\lambda x_k, x_j \right)_{\widehat H} $$
$$ = \left( U^{-1} \int_{\mathbb{R}_+} \frac{1}{\lambda - z} d\widehat E_\lambda x_k, U^{-1} x_j \right)_{L^2(\widehat M)} $$
$$ = \left( \int_{\mathbb{R}_+} \frac{1}{\lambda - z} d U^{-1} \widehat E_\lambda U \vec e_k, \vec e_j \right)_{L^2(\widehat M)} =
\left( \int_{\mathbb{R}_+} \frac{1}{\lambda - z} d E_{\lambda} \vec e_k, \vec e_j \right)_{L^2(\widehat M)} $$
\begin{equation}
\label{f3_10}
= \int_{\mathbb{R}_+} \frac{1}{\lambda - z} d(E_{\lambda} \vec e_k, \vec e_j)_{L^2(\widehat M)},\qquad 0\leq k,j\leq N-1.
\end{equation}
Using~(\ref{f3_2}) we can write
$$ (E_{\lambda} \vec e_k, \vec e_j)_{L^2(\widehat M)} = \widehat m_{k,j}(\lambda), $$
and therefore
\begin{equation}
\label{f3_11}
\int_{\mathbb{R}_+} \frac{1}{\lambda - z} d( P^{\widehat H}_H \widehat E_\lambda x_k, x_j)_H =
\int_{\mathbb{R}_+} \frac{1}{\lambda - z} d\widehat m_{k,j}(\lambda),\qquad 0\leq k,j\leq N-1.
\end{equation}
By the Stieltjes-Perron inversion  formula (see, e.g.,~\cite{c_14000_Akh}) we conclude that
\begin{equation}
\label{f3_12}
\widehat m_{k,j} (\lambda) = ( P^{\widehat H}_H \widehat E_\lambda x_k, x_j)_H.
\end{equation}

\begin{prop}
\label{p3_1}
Let the matrix Stieltjes moment problem~(\ref{f1_1}) be given and conditions~(\ref{f1_4}) hold.
Let $A$ be a non-negative Hermitian operator which is defined by~(\ref{f2_12}).
The deficiency index of $A$ is equal to $(n,n)$, $0\leq n\leq N$.
\end{prop}
{\bf Proof. }
Choose an arbitrary $u\in L$, $u = \sum_{k=0}^\infty c_k x_k$, $c_k\in \mathbb{C}$. Suppose that
$c_k = 0$, $k\geq N+R+1$, for some $R\in \mathbb{Z}_+$. Consider the following system of linear equations:
\begin{equation}
\label{f3_13}
-z d_k = c_k,\qquad  k=0,1,...,N-1;
\end{equation}
\begin{equation}
\label{f3_14}
d_{k-N} - z d_k = c_k,\qquad  k=N,N+1,N+2,...;
\end{equation}
where $\{ d_k \}_{k\in \mathbb{Z}_+}$ are unknown complex numbers, $z\in \mathbb{C}\backslash \mathbb{R}$ is a
fixed parameter.
Set
$$ d_k = 0,\qquad k\geq R+1; $$
\begin{equation}
\label{f3_15}
d_{j} = c_{N+j} + z d_{N+j},\qquad j=R,R-1,R-2,...,0.
\end{equation}
For such defined numbers $\{ d_k \}_{k\in\mathbb{Z}_+}$, all equations in~(\ref{f3_14}) are satisfied.
But equations~(\ref{f3_14}) are not necessarily satisfied. Set
$$ v = \sum_{k=0}^\infty d_k x_k,\ v\in L. $$
Notice that
$$ (A-zE_H) v = \sum_{k=0}^\infty (d_{k-N} - z d_k) x_k, $$
where $d_{-1}=d_{-2}=...=d_{-N}=0$.
By the construction of $d_k$ we have
$$ (A-zE_H) v - u = \sum_{k=0}^\infty (d_{k-N} - z d_k - c_k) x_k =
\sum_{k=0}^{N-1} (-zd_k - c_k) x_k; $$
\begin{equation}
\label{f3_16}
u = (A-zE_H) v + \sum_{k=0}^{N-1} (zd_k + c_k) x_k,\qquad u\in L.
\end{equation}
Set
$$ H_z := \overline{(A-zE_H) L} = (\overline{A} - zE_H) D(\overline{A}), $$
and
\begin{equation}
\label{f3_17}
y_k := x_k - P^H_{H_z} x_k,\qquad k=0,1,...,N-1.
\end{equation}
Set
$$ H_0 := \mathop{\rm span}\nolimits\{ y_k \}_{k=0}^{N-1}. $$
Notice that the dimension of $H_0$ is less or equal to $N$, and $H_0\perp H_z$.
From~(\ref{f3_16}) it follows that $u\in L$ can be represented in the following form:
\begin{equation}
\label{f3_18}
u = u_1 + u_2,\qquad u_1\in H_z,\quad u_2\in H_0.
\end{equation}
Therefore we get $L\subseteq H_z\oplus H_0$; $H\subseteq H_z\oplus H_0$, and finally
$H=H_z\oplus H_0$. Thus, $H_0$ is the corresponding defect subspace.
So, the defect numbers of $A$ are less or equal to $N$. Since the operator $A$ is non-negative,
they are equal.
$\Box$

\begin{thm}
\label{t3_1}
Let a matrix Stieltjes moment problem~(\ref{f1_1}) be given and
conditions~(\ref{f1_4}) hold. Let an operator $A$ be constructed for the
moment problem as in~(\ref{f2_12}).
All solutions of the moment problem have the following form
\begin{equation}
\label{f3_19}
M(\lambda) = (m_{k,j} (\lambda))_{k,j=0}^{N-1},\quad
m_{k,j} (\lambda) = ( \mathbf E_\lambda x_k, x_j)_H,
\end{equation}
where $\mathbf E_\lambda$ is a $\Pi$-spectral function of the operator $A$.
Moreover, the correspondence between all $\Pi$-spectral functions of $A$ and all solutions
of the moment problem is one-to-one.
\end{thm}
{\bf Proof. }
It remains to prove that different $\Pi$-spectral functions of the operator $A$ produce different
solutions of the moment problem~(\ref{f1_1}).
Suppose to the contrary that two different $\Pi$-spectral functions produce the same solution of
the moment problem. That means that
there exist two non-negative self-adjoint extensions
$A_j\supseteq A$, in Hilbert spaces $H_j\supseteq H$, such that
\begin{equation}
\label{f3_20}
P_{H}^{H_1} E_{1,\lambda} \not= P_{H}^{H_2} E_{2,\lambda},
\end{equation}
\begin{equation}
\label{f3_21}
(P_{H}^{H_1} E_{1,\lambda} x_k,x_j)_H = (P_{H}^{H_2} E_{2,\lambda} x_k,x_j)_H,\qquad 0\leq k,j\leq N-1,\quad
\lambda\in\mathbb{R}_+,
\end{equation}
where $\{ E_{n,\lambda} \}_{\lambda\in\mathbb{R}_+}$ are orthogonal left-continuous resolutions of unity of
operators $A_n$, $n=1,2$.
Set $L_N := \mathop{\rm Lin}\nolimits\{ x_k \}_{k=0,N-1}$. By linearity we get
\begin{equation}
\label{f3_22}
(P_{H}^{H_1} E_{1,\lambda} x,y)_H = (P_{H}^{H_2} E_{2,\lambda} x,y)_H,\qquad x,y\in L_N,\quad \lambda\in\mathbb{R}_+.
\end{equation}
Denote by $R_{n,\lambda}$ the resolvent of $A_n$, and set $\mathbf R_{n,\lambda} := P_{H}^{H_n} R_{n,\lambda}$, $n=1,2$.
From~(\ref{f3_22}),(\ref{f3_1}) it follows that
\begin{equation}
\label{f3_23}
(\mathbf R_{1,z} x,y)_H = (\mathbf R_{2,z} x,y)_H,\qquad x,y\in L_N,\quad z\in \mathbb{C}\backslash \mathbb{R}.
\end{equation}
Choose an arbitrary $z\in\mathbb{C}\backslash\mathbb{R}$ and consider the space $H_z$ defined as above.
Since
$$ R_{j,z} (A-zE_H) x = (A_j - z E_{H_j} )^{-1} (A_j - z E_{H_j}) x = x,\qquad x\in L=D(A),$$
we get
\begin{equation}
\label{f3_24}
R_{1,z} u = R_{2,z} u \in H,\qquad u\in H_z;
\end{equation}
\begin{equation}
\label{f3_25}
\mathbf R_{1,z} u = \mathbf R_{2,z} u,\qquad u\in H_z,\ z\in\mathbb{C}\backslash\mathbb{R}.
\end{equation}
We can write
$$ (\mathbf R_{n,z} x, u)_H = (R_{n,z} x, u)_{H_n} = ( x, R_{n,\overline{z}}u)_{H_n} =
( x, \mathbf R_{n,\overline{z}} u)_H, $$
\begin{equation}
\label{f3_26}
x\in L_N,\ u\in H_{\overline z},\ n=1,2,
\end{equation}
and therefore we get
\begin{equation}
\label{f3_27}
(\mathbf R_{1,z} x,u)_H = (\mathbf R_{2,z} x,u)_H,\qquad x\in L_N,\ u\in H_{\overline z}.
\end{equation}
By~(\ref{f3_16}) an arbitrary element $y\in L$ can be represented as $y=y_{ \overline{z} } + y'$,
$y_{ \overline{z} }\in H_{ \overline{z} }$, $y'\in L_N$.
Using~(\ref{f3_23}) and~(\ref{f3_25})  we get
$$ (\mathbf R_{1,z} x,y)_H = (\mathbf R_{1,z} x, y_{ \overline{z} } + y')_H $$
$$ = (\mathbf R_{2,z} x, y_{ \overline{z} } + y')_H = (\mathbf R_{2,z} x,y)_H,\qquad x\in L_N,\ y\in L. $$
Since $\overline{L}=H$, we obtain
\begin{equation}
\label{f3_28}
\mathbf R_{1,z} x = \mathbf R_{2,z} x,\qquad x\in L_N,\ z\in\mathbb{C}\backslash\mathbb{R}.
\end{equation}
For an arbitrary $x\in L$, $x=x_z + x'$, $x_z\in H_z$, $x'\in L_N$, using
relations~(\ref{f3_25}),(\ref{f3_28}) we obtain
\begin{equation}
\label{f3_29}
\mathbf R_{1,z} x = \mathbf R_{1,z} (x_z + x') =
\mathbf R_{2,z} (x_z + x') = \mathbf R_{2,z} x,\qquad x\in L,\ z\in\mathbb{C}\backslash\mathbb{R},
\end{equation}
and
\begin{equation}
\label{f3_30}
\mathbf R_{1,z} x = \mathbf R_{2,z} x,\qquad x\in H,\ z\in\mathbb{C}\backslash\mathbb{R}.
\end{equation}
By~(\ref{f3_1}) that means that the $\Pi$-spectral functions coincide and we obtain a
contradiction.
$\Box$

We shall recall some basic definitions and facts from~\cite{c_9000_D_M}.
Let $A$ be a closed Hermitian operator in a Hilbert space $H$, $\overline{D(A)} = H$.
\begin{dfn}
\label{d3_1}
A collection $\{ \mathcal{H}, \Gamma_1, \Gamma_2 \}$ in which $\mathcal{H}$ is a Hilbert space,
$\Gamma_1, \Gamma_2 \in [D(A^*),\mathcal{H}]$, is called {\bf a space of boundary values (SBV)} for
$A^*$, if

\noindent
(1) $(A^* f,g)_H - (f,A^* g)_H = (\Gamma_1 f,\Gamma_2 g)_{\mathcal{H}} - (\Gamma_2 f, \Gamma_1 g)_{\mathcal{H}}$,
$\forall f,g\in D(A^*)$;

\noindent
(2) the mapping $\Gamma: f\rightarrow \{ \Gamma_1 f,\Gamma_2 f \}$ from $D(A^*)$ to $\mathcal{H}\oplus
\mathcal{H}$ is surjective.
\end{dfn}
Naturally associated with each SBV are self-adjoint operators $\widetilde A_1,\widetilde A_2\ (\subset A^*)$ with
$$ D(\widetilde A_1) = \mathop{\rm ker}\nolimits \Gamma_1,\ D(\widetilde A_2) = \mathop{\rm ker}\nolimits \Gamma_2. $$
The operator $\Gamma_2$ restricted to the defect subspace $N_z = \mathop{\rm ker}\nolimits(A^* - zE_H)$,
$z\in \rho(\widetilde A_2)$, is fully invertible. For
$\forall z\in \rho(\widetilde A_2)$ set
\begin{equation}
\label{f3_31}
\gamma(z) = \left( \Gamma_2|_{N_z} \right)^{-1} \in [\mathcal{H},N_z].
\end{equation}
\begin{dfn}
\label{d3_2}
The operator-valued function $M(z)$ defined for $z\in\rho(\widetilde A_2)$ by
\begin{equation}
\label{f3_32}
M(z) \Gamma_2 f_z = \Gamma_1 f_z,\qquad f_z\in N_z,
\end{equation}
is called {\bf a Weyl function} of the  operator $A$, corresponding to SBV $\{ \mathcal{H}, \Gamma_1, \Gamma_2 \}$.
\end{dfn}
The Weyl function can be also obtained from the equality:
\begin{equation}
\label{f3_33}
M(z) = \Gamma_1 \gamma(z),\qquad z\in \rho(\widetilde A_2).
\end{equation}
For an arbitrary operator $\widetilde A = \widetilde A^* \subset A^*$ there exist a SBV with (\cite{c_15000_D_M})
\begin{equation}
\label{f3_34}
D(\widetilde A_2) = \mathop{\rm ker}\nolimits\Gamma_2 = D(\widetilde A).
\end{equation}
(There even exist a family of such SBV).

An extension $\widehat A$ of $A$ is called {\bf proper} if $A\subset\widehat A\subset A^*$ and
$(\widehat A^*)^* = \widehat A$. Two proper extensions $\widehat A_1$ and $\widehat A_2$ are
{\bf  disjoint} if $D(\widehat A_1)\cap D(\widehat A_2) = D(A)$ and
{\bf transversals} if they are disjoint and $D(\widehat A_1) + D(\widehat A_2) = D(A^*)$.

Suppose that the operator $A$ is non-negative, $A\geq 0$. In this case there exist
two non-negative self-adjoint extensions of $A$ in $H$, Friedrich's extension $A_\mu$ and Krein's extension $A_M$, such that
for an arbitrary non-negative self-adjoint extension $\widehat A$ of $A$ in $H$ it holds:
\begin{equation}
\label{f3_35}
(A_\mu + xE_H)^{-1} \leq (\widehat A + xE_H)^{-1} \leq (A_M + xE_H)^{-1},\qquad x\in \mathbb{R}_+.
\end{equation}
Recall some definitions and facts from~\cite{c_8000_Kr_Ovch},\cite{c_13000_Kr}.
For the non-negative operator $A$ we put into correspondence the following operator:
\begin{equation}
\label{f3_36}
T = (E_H - A)(E_H + A)^{-1} = -E_H + 2(E_H + A)^{-1},\qquad D(T) = (A+E_H)D(A).
\end{equation}
The operator $T$ is a Hermitian contraction (i.e. $\| T \| \leq 1$). Its domain is not dense in $H$
if $A$ is not self-adjoint. The defect subspace $H\ominus D(T) = N_{-1}$ and its dimension is equal
to the defect number $n(A)$ of $A$.
The inverse transformation to~(\ref{f3_36}) is given by
\begin{equation}
\label{f3_37}
A = (E_H - T)(E_H + T)^{-1} = -E_H + 2(E_H + T)^{-1},\qquad D(A) = (T+E_H)D(T).
\end{equation}
Relations~(\ref{f3_36}),(\ref{f3_37}) (with $\widehat T,\widehat A$ instead of $T,A$) also
establish a bijective correspondence between
self-adjoint contractive extensions $\widehat T\supseteq T$ in $H$ and self-adjoint non-negative
extensions $\widehat A\supseteq A$ in $H$ (\cite[p.451]{c_13000_Kr}).

\noindent
Consider an arbitrary Hilbert space $\widehat H \supseteq H$. It is not hard to see that
relations~(\ref{f3_36}),(\ref{f3_37}) (with $\widehat T,\widehat A$ instead of $T,A$) establish
a bijective correspondence between
self-adjoint contractive extensions $\widehat T\supseteq T$ in $\widehat H$ and self-adjoint non-negative
extensions $\widehat A\supseteq A$ in $\widehat H$, as well.

There exist extremal self-adjoint contractive extensions of $T$ in $H$ such that for
an arbitrary self-adjoint contractive extension $\widetilde T\supseteq T$
in $H$ it holds
\begin{equation}
\label{f3_38}
T_\mu \leq \widetilde T \leq T_M.
\end{equation}
Notice that
\begin{equation}
\label{f3_39}
A_\mu = -E_H + 2(E_H + T_\mu)^{-1},\quad  A_M = -E_H + 2(E_H + T_M)^{-1}.
\end{equation}
Set
\begin{equation}
\label{f3_40}
C = T_M - T_\mu.
\end{equation}
Consider the following subspace:
\begin{equation}
\label{f3_41}
\Upsilon = \mathop{\rm ker}\nolimits \left( C|_{ N_{-1} } \right).
\end{equation}
\begin{dfn}
\label{d3_3}
Let a closed non-negative Hermitian operator $A$ be  given. For the operator $A$ it takes place
{\bf a completely indeterminate case} if $\Upsilon = \{ 0 \}$.
\end{dfn}
By Theorem~1.4 in~\cite{c_16000_KO}, on the set $\{ x\in H:\ T_\mu x = T_M x \} = \mathop{\rm ker}\nolimits C$,
all self-adjoint contractive extensions in a Hilbert space $\widetilde H\supseteq H$ coincide.
Thus, all such extensions are extensions of the operator $T_{ext}$:
\begin{equation}
\label{f3_42}
T_{ext} x = \left\{\begin{array}{cc} Tx, & x\in D(T)\\
T_\mu x = T_M x, & x\in \mathop{\rm ker}\nolimits C \end{array}\right..
\end{equation}
Introduce the following operator:
\begin{equation}
\label{f3_43}
A_{ext} = -E_H + 2(E_H + T_{ext})^{-1} \supseteq A.
\end{equation}
Thus, {\it the set of all non-negative self-adjoint extensions of $A$ coincides with the set of
all non-negative self-adjoint extensions of $A_{ext}$.}
Since $T_{ext,\mu} = T_\mu$ and $T_{ext,M} = T_M$, for $A_{ext}$ it takes place
the completely indeterminate case.

\begin{prop}
\label{p3_2}
Let $A$ be a closed non-negative Hermitian operator with finite defect numbers
and for $A$ it takes place the completely indeterminate case.
Then extensions $A_\mu$ and $A_M$ given by~(\ref{f3_39}) are transversal.
\end{prop}
{\bf Proof. }
Notice that
\begin{equation}
\label{f3_44}
D(A_M) \cap D(A_\mu) = D(A).
\end{equation}
In fact, suppose that there exists $y\in D(A_M) \cap D(A_\mu)$, $y\notin D(A)$.
Since $A_M\subset A^*$ and $A_\mu\subset A^*$ we have $A_M y = A_\mu y$.
Set
$$ g := (A_M + E_H)y = (A_\mu + E_H)y . $$
Then
$$ T_M g = -g + 2(E_H + A_M)^{-1}g = -g + 2y, $$
$$ T_\mu g = -g + 2(E_H + A_\mu)^{-1}g = -g + 2y, $$
and therefore $Cg = (T_M - T_\mu)g = 0$. Since $y\notin D(A)$, then $g\in N_{-1}$.
We obtained a contradiction, since for $A$ it takes place the completely indeterminate case.

Introduce the following sets:
\begin{equation}
\label{f3_45}
D_M := (A_M + E_H)^{-1} N_{-1},\quad D_\mu := (A_\mu + E_H)^{-1} N_{-1}.
\end{equation}
Since $D(A_M) = (A_M + E_H)^{-1} D(T_M)$, $D(A_\mu) = (A_\mu + E_H)^{-1} D(T_\mu)$, we have
\begin{equation}
\label{f3_46}
D_M \subset D(A_M),\  D_\mu \subset D(A_\mu),
\end{equation}
and
\begin{equation}
\label{f3_47}
D_M \cap D(A) = \{ 0 \},\quad  D_\mu \cap D(A) = \{ 0 \},
\end{equation}
By~(\ref{f3_44}),(\ref{f3_46}) and~(\ref{f3_47}) we obtain that
\begin{equation}
\label{f3_48}
D_M \cap D_\mu = \{ 0 \}.
\end{equation}
Set
\begin{equation}
\label{f3_49}
D := D_M \dotplus D_\mu.
\end{equation}
By~(\ref{f3_45}) we obtain that the sets $D_M$ and $D_\mu$ have the linear dimension $n(A)$.
Elementary arguments show that $D$ has the linear dimension $2 n(A)$.
Since $D(A_\mu)\subset D(A^*)$, $D(A_M)\subset D(A^*)$, we can write
\begin{equation}
\label{f3_50}
D(A) \dotplus D_M \dotplus D_\mu \subseteq D(A^*) = D(A) \dotplus N_z \dotplus N_{\overline{z}},
\end{equation}
where $z\in \mathbb{C}\backslash \mathbb{R}$.

Let
$$ g_1,g_2,...,g_{2n(A)}, $$
be $2n(A)$ linearly independent elements from $D$. Let
\begin{equation}
\label{f3_52}
g_j = g_{A,j} + g_{z,j} + g_{\overline{z},j},\qquad 1\leq j\leq 2n(A),
\end{equation}
where $g_{A,j}\in D(A)$, $g_{z,j}\in N_z$, $g_{\overline{z},j}\in N_{\overline{z}}$.
Set
\begin{equation}
\label{f3_53}
\widehat g_j := g_j - g_{A,j},\qquad 1\leq j\leq 2n(A).
\end{equation}
If for some $\alpha_j\in \mathbb{C}$, $1\leq j\leq 2n(A)$, we have
$$ 0 = \sum_{j=1}^{ 2n(A) } \alpha_j \widehat g_j = \sum_{j=1}^{ 2n(A) } \alpha_j g_j -
\sum_{j=1}^{ 2n(A) } \alpha_j  g_{A,j}, $$
then
$$ \sum_{j=1}^{ 2n(A) } \alpha_j g_j = 0, $$
and $\alpha_j = 0$, $1\leq j\leq 2n(A)$.
Therefore elements $\widehat g_j$, $1\leq j\leq 2n(A)$ are linearly independent.
Thus, they form a linear basis in a finite-dimensional subspace $N_z \dotplus N_{\overline{z}}$.
Then
\begin{equation}
\label{f3_54}
N_z \dotplus N_{\overline{z}} \subseteq D,
\end{equation}
\begin{equation}
\label{f3_55}
D(A^*) = D(A) \dotplus N_z \dotplus N_{\overline{z}} \subseteq D(A) \dotplus D = D_L.
\end{equation}
So, we get the equality
\begin{equation}
\label{f3_56}
D(A) \dotplus D_M \dotplus D_\mu = D(A^*).
\end{equation}

Since $D(A) + D_M \subseteq D(A_M)$, $D_\mu \subseteq D(A_\mu)$,
we get
$$ D(A^*) = D(A)+D_M+D_\mu \subseteq D(A_M) + D(A_\mu). $$
Since $D(A_M) + D(A_\mu)\subseteq D(A^*)$, we get
\begin{equation}
\label{f3_57}
D(A^*) = D(A_M) + D(A_\mu).
\end{equation}
From~(\ref{f3_44}),(\ref{f3_57}) it follows the statement of the Proposition.
$\Box$

We shall use the following classes of functions~\cite{c_9000_D_M}.
Let $\mathcal{H}$ be a Hilbert space.
Denote by $R_{\mathcal{H}}$ the class of operator-valued functions $F(z)=F^*(\overline{z})$
holomorphic in $\mathbb{C}\backslash \mathbb{R}$ with values (for $z\in \mathbb{C}_+$)
in the set of maximal dissipative operators in $\mathfrak{C}(\mathcal{H})$.
Completing the class $R_{\mathcal{H}}$ by ideal elements we get the class $\widetilde R_{\mathcal{H}}$.
Thus, $\widetilde R_{\mathcal{H}}$ is a collection of functions holomorphic in
$\mathbb{C}\backslash \mathbb{R}$ with values (for $z\in \mathbb{C}_+$)
in the set of maximal dissipative linear relations $\theta(z)=\theta^*(\overline{z})$ in $\mathcal{H}$.
The indeterminate part of the relation $\theta(z)$ does not depend on $z$ and the relation $\theta(z)$
admits the representation
\begin{equation}
\label{f3_58}
\theta(z) = \{ < h_1, F_1(z) h_1 + h_2 >:\ h_1\in D(F_1(z)),\ h_2\in \mathcal{H}_2 \},
\end{equation}
where $\mathcal{H} = \mathcal{H}_1 \oplus \mathcal{H}_2$, $F_1(z) \in R_{\mathcal{H}_1}$.
\begin{dfn} \cite{c_9000_D_M}
\label{d3_4}
An operator-valued function $F(z)\in R_{\mathcal{H}}$ belongs to the class
$S_{\mathcal{H}}^{-0} (-\infty,0)$ if $\forall n\in \mathbb{N}$, $\forall z_j\in \mathbb{C}_+$,
$h_j\in D(F(z_j))$, $\xi_j\in \mathbb{C}$, holds
\begin{equation}
\label{f3_59}
\sum_{i,j=1}^n \frac{ (z_i^{-1} F(z_i) h_i, h_j) - (h_i, z_j^{-1} F(z_j) h_j) }
{ z_i - \overline{z_j} } \xi_i \overline{\xi_j} \geq 0.
\end{equation}
Completing the class $S_{\mathcal{H}}^{-0} (-\infty,0)$ with ideal elements~(\ref{f3_58})
we obtain the class $\widetilde S_{\mathcal{H}}^{-0} (-\infty,0)$.
\end{dfn}
From Theorem~9 in~\cite[p.46]{c_9000_D_M} taking into account Proposition~\ref{p3_2} we have the
following conclusion (see also Remark~17 in~\cite[p.49]{c_9000_D_M}):
\begin{thm}
\label{t3_2}
Let $A$ be a closed non-negative Hermitian operator in a Hilbert space $H$ and for $A$ it takes place
the completely indeterminate case.
Let $\{ \mathcal{H}, \Gamma_1, \Gamma_2 \}$ be an arbitrary SBV for $A$ such that
$\widetilde A_2 = A_\mu$ and $M(z)$ be the corresponding Weyl function. Then the formula
\begin{equation}
\label{f3_60}
\mathbf{R}_z = (A_\mu - zE_H)^{-1} - \gamma(z) (\tau(z)+M(z)-M(0))^{-1} \gamma^*(\overline{z}),\quad
z\in \mathbb{C}\backslash \mathbb{R},
\end{equation}
establishes a bijective correspondence between $\mathbf{R}_z\in \Omega^0(-\infty,0)(A)$ and
$\tau\in \widetilde S_{\mathcal{H}}^{-0} (-\infty,0)$.
The function $\tau(z)\equiv \tau = \tau^*$ in~(\ref{f3_60}) corresponds to the canonical $\Pi$-resolvents
and only to them.
\end{thm}
Now we can state our main result.
\begin{thm}
\label{t3_3}
Let a matrix Stieltjes moment problem~(\ref{f1_1}) be given and
conditions~(\ref{f1_4}) hold. Let an operator $A$ be the closure of the operator constructed for the
moment problem in~(\ref{f2_12}).
Then the following statements are true:

1) The moment problem~(\ref{f1_1}) is determinate if and only if
Friedrich's extension $A_\mu$ and Krein's extension $A_M$ coincide: $A_\mu = A_M$.
In this case the unique solution of the moment problem is generated by the orthogonal spectral
function $\mathbf E_\lambda$ of $A_\mu$ by formula~(\ref{f3_19});

2) If $A_\mu\not= A_M$, define the extended operator $A_{ext}$ for $A$ as in~(\ref{f3_43}).
Let $\{ \mathcal{H}, \Gamma_1, \Gamma_2 \}$ be an arbitrary SBV for $A_{ext}$ such that
$\widetilde A_2 = (A_{ext})\mu$ and $M(z)$ be the corresponding Weyl function.
All solutions of the moment problem~(\ref{f1_1}) have the following form:
\begin{equation}
\label{f3_61}
M(\lambda) = (m_{k,j} (\lambda))_{k,j=0}^{N-1},
\end{equation}
where
$$ \int_{\mathbb{R}_+} \frac{dm_{k,j} (\lambda)}{\lambda - z} =
\left( (A_\mu - zE_H)^{-1} x_k,x_j \right)_H $$
\begin{equation}
\label{f3_62}
- \left( \gamma(z) (\tau(z)+M(z)-M(0))^{-1} \gamma^*(\overline{z}) x_k,
x_j \right)_H,\quad z\in \mathbb{C}\backslash \mathbb{R},
\end{equation}
where $\tau\in \widetilde S_{\mathcal{H}}^{-0} (-\infty,0)$.
Moreover, the correspondence between all $\tau\in \widetilde S_{\mathcal{H}}^{-0} (-\infty,0)$ and all solutions
of the moment problem~(\ref{f1_1}) is one-to-one.
\end{thm}
{\bf Proof. } The statements of the Theorem follow directly from Theorems~\ref{t3_1} and~\ref{t3_2}.
$\Box$

\begin{center}
\large\bf
The matrix Stieltjes moment problem: a description of all solutions.
\end{center}
\begin{center}
\bf S.M. Zagorodnyuk
\end{center}

We describe all solutions of the matrix Stieltjes moment problem in the general case (no conditions
besides solvability are assumed). We use Krein's formula for the generalized $\Pi$-resolvents
of positive Hermitian operators in the form of V.A~Derkach and M.M.~Malamud.

MSC: 44A60; Secondary 30E05
Key words:  moment problem, positive definite kernel, spectral function.

\end{document}